\documentclass{article}[12pt]
\usepackage{amssymb}
\usepackage{eurosym}
\usepackage{epsfig}
\usepackage{amsmath}
\usepackage{amsfonts}
\usepackage{pgf,tikz}
\usepackage{enumerate}
\usepackage{cite}
\usepackage{xcolor, colortbl}

\usepackage[colorlinks=true, linkcolor=black, citecolor=black, urlcolor=black]{hyperref}

\newtheorem{theorem}{\sc Theorem}[section]
\newtheorem{proposition}{\sc Proposition}[section]
\newtheorem{lemma}{\sc Lemma}[section]

\newtheorem{remark}{\sc Remark}[section]
\newtheorem{corollary}{\sc Corollary}[section]
\newtheorem{example}{\sc Example}[section]

\def\qed{\hbox to 0pt{}\hfill$\rlap{$\sqcap$}\sqcup$\medbreak}

\title{Component-wise Krasnosel'ski\u{\i} type fixed point theorem in product spaces and applications}

\author{Laura Mª Fern\'andez--Pardo and Jorge Rodr\'iguez--L\'opez} 

\date{}
\begin{document}
 \maketitle

\begin{center}  {\small CITMAga \& Departamento de Estat\'{\i}stica, An\'alise Matem\'atica e Optimizaci\'on, \\ Universidade de Santiago de Compostela, \\ 15782, Facultade de Matem\'aticas, Campus Vida, Santiago, Spain.\\  Email: lauramaria.fernandez.pardo@rai.usc.es; jorgerodriguez.lopez@usc.es}
\end{center}

\medbreak

\noindent {\it Abstract.} We present a version of Krasnosel’ski\u{\i} fixed point theorem for operators acting on Cartesian products of normed linear spaces, under cone-compression and cone-expansion conditions of norm type. Our approach, based on the fixed point index theory in cones, guarantees the existence of a coexistence fixed point—that is, one with nontrivial components. As an application, we prove the existence of periodic solutions with strictly positive components for a system of second-order differential equations. In particular, we address cases involving singular nonlinearities and hybrid terms, characterized by sublinear behavior in one component and superlinear behavior in the other. 

\medbreak


\medbreak


\section{Introduction}

Fixed point theorems for completely continuous maps in normed linear spaces play a crucial role in the study of nonlinear boundary value problems, see for instance \cite{amann,guolak}. In particular, the cone compression-expansion fixed point theorems are commonly employed in order to prove the existence of positive non-trivial solutions. In this paper, we focus on Krasnosel'ski\u{\i} fixed point theorem with norm type compression/expansion conditions, which localizes the solution in a conical shell.

In the case of systems of operator equations, this type of localization does not prevent the obtained solution from being semi-trivial, in the sense that some of its components may be zero, as already pointed out in \cite{PrecupFPT,PrecupSDC}. Therefore, this result does not allow for the direct determination of \textit{coexistence fixed points}, as coined by Lan \cite{Lan1}. This observation motivated Precup \cite{PrecupFPT,PrecupSDC} to introduce vectorial extensions of Krasnosel’ski\u{\i} fixed point theorem under compression and expansion conditions. Specifically, one of his approaches imposes constraints based on the partial order induced by the cone, while the other is formulated in terms of homotopy-type conditions. Alternative results concerning fixed points which are non-trivial in all the components can be found in \cite{ima,Lan1,JRL}. Nevertheless, as far as we are aware, there is no such a vectorial version of Krasnosel'ski\u{\i} theorem with norm type cone-compression and cone-expansion conditions. Our main abstract result in this paper, which is the following fixed point theorem in Cartesian products, fills this gap in the literature.         

\begin{theorem}\label{main_th}
	Let $K_1$ and $K_2$ be two cones in the normed spaces $X_1$ and $X_2$, respectively. We denote $K:=K_1\times K_2$ the corresponding cone in $X=X_1\times X_2$, and for $r,R\in\mathbb{R}^2$ with $0<r_i<R_i$ ($i=1,2$) 
	 \[\overline{K}_{r,R}:=\left\{ x=(x_1,x_2)\in K: r_i\leq\left\|x_i\right\|\leq R_i \text{ for } i=1,2\right\}.\]
	Assume that $\mathcal{T}=(\mathcal{T}_1,\mathcal{T}_2):\overline{K}_{r,R}\rightarrow K$ is a compact map and for each $i\in\{1,2\}$ either of the following conditions holds in $\overline{K}_{r,R}$:
	\begin{enumerate}
		\item[(a)] $\left\|\mathcal{T}_i x\right\|\geq \left\|x_i\right\|$ if $\left\|x_i\right\|=r_i$ and  $\left\|\mathcal{T}_i x\right\|\leq \left\|x_i\right\|$ if $\left\|x_i\right\|=R_i$; or
		\item[(b)] $\left\|\mathcal{T}_i x\right\|\leq \left\|x_i\right\|$ if $\left\|x_i\right\|=r_i$ and  $\left\|\mathcal{T}_i x\right\|\geq \left\|x_i\right\|$ if $\left\|x_i\right\|=R_i$. 
	\end{enumerate}
	Then $\mathcal{T}$ has at least one fixed point $x=(x_1,x_2)\in K$ such that $r_i\leq \left\|x_i\right\|\leq R_i$ ($i=1,2$).
\end{theorem}

In the sequel, following the usual terminology, it is said that the component $\mathcal{T}_i$ of the operator $\mathcal{T}$ is \textit{compressive} if $(a)$ holds, whereas $\mathcal{T}_i$ is said to be \textit{expansive} in case that condition $(b)$ is satisfied.  
When compared with the classical Krasnosel'ski\u{\i} fixed point theorem in cones, Theorem~\ref{main_th} has a twofold interest: 
\begin{enumerate}
\item[$(i)$] it ensures the existence of a coexistence fixed point; 
\item[$(ii)$] different behaviors (compressive or expansive) are allowed for each of the components of the operator.
\end{enumerate}
Note that, in particular, this second advantage allows to deal with expansive-compressive maps, which have not been extensively studied in the literature. We refer the interested reader to the manuscripts by Mawhin \cite{Maw_CE}, Pireddu-Zanolin \cite{pz} and the references therein.

For the proof of Theorem \ref{main_th}, we consider separately the case in which both components, $\mathcal{T}_1$ and $\mathcal{T}_2$, are compressive. In this setting, the reasoning is based on adequate computations of the Leray-Schauder fixed point index for compact maps in cones. It is noteworthy that obtaining the value of the fixed-point index facilitates the derivation of multiplicity results, analogous to those established by Amann \cite{amann}. Finally, the remaining possibilities for the behaviors of $\mathcal{T}_1$ and $\mathcal{T}_2$ can be reduced to the compressive case.  

As an application, we are concerned with the existence of $T$-periodic solutions with strictly positive components for a system of second-order equations of the form
\begin{equation}\label{eq_p}
	\left\{\begin{array}{l} x''+a_1(t)\,x=f_1(t,x,y), \\[6pt] y''+a_2(t)\,y=f_2(t,x,y), \end{array} \right.
\end{equation}
where $a_1$ and $a_2$ guarantee that the problem is non-resonant and the corresponding Green's functions has constant sign, and the nonlinearities $f_1$ and $f_2$ are continuous. Our existence results provide a way to adapt those due to Torres \cite{Torres}, in the scalar case, to the setting of systems. As a novelty, we highlight that the expansive/compressive behavior mentioned in $(ii)$, when applied to system \eqref{eq_p}, allows to deal with right-hand sides which are sublinear in one component and superlinear in the other one, see Example~\ref{ex_period_subsup}. In addition, the case of \textit{singular} nonlinearities is also considered, complementing previous existence criteria available in the literature, see \cite{W,WCS,ft} and the references therein. In particular, we deal with singularities which are \textit{repulsive} in one component and \textit{attractive} in the other, in the sense specified in Remark \ref{rmk_a_r}. 

The paper is organized as follows. In Section~\ref{sec2}, we recall some basic properties of the Leray-Schauder fixed point index in cones and, moreover, we show a specific computation of the index for operators defined in Cartesian products, which is crucial in the proof of our main results. In Section~\ref{sec_main}, we prove Theorem~\ref{main_th} and, furthermore, we establish multiplicity results. Finally, Section~\ref{sec_app} is devoted to the existence of positive periodic solutions for a system of second-order equations.     

\section{Fixed point index computations}\label{sec2}
Let us recall some concepts and properties related to the fixed point index in cones. A closed convex subset $K$ of a normed linear space $X$ is said to be a \textit{cone} if $\lambda\,u\in K$ for every $u\in K$ and for all $\lambda\geq 0$, and $K\cap (-K)=\{0\}$.
If $U$ is a relatively open bounded subset of a cone $K$ of a normed space $X$ and $N:\overline{U}\rightarrow K$ is a compact map without fixed points on the boundary of $U$ (denoted by $\partial_K\,U$), the \textit{fixed point index} of $N$ on $U$ with respect to the cone $K$, $i_K(N,U)$, is well-defined. For details, we refer the reader to \cite{amann,guolak}. 


\begin{proposition}\label{prop_index}
	Let $K$ be a cone of a normed space, $U\subset K$ be a bounded relatively open set and $N:\overline{U}\rightarrow K$ be a compact map such that $N$ has no fixed points on the boundary of $U$. Then the fixed point index of $N$ on the set $U$ with respect to $K$, $i_{K}(N,U)$, has the following properties:
	\begin{enumerate}
		\item (Additivity) Let $U$ be the disjoint union of two open sets $U_1$ and $U_2$. If $0\not\in(I-N)(\overline{U}\setminus(U_1\cup U_2))$, then \[i_{K}(N,U)=i_{K}(N,U_1)+i_{K}(N,U_2).\]
		\item (Existence) If $i_{K}(N,U)\neq 0$, then there exists $x\in U$ such that $Nx=x$.
		\item (Homotopy invariance) If $H:\overline{U}\times[0,1]\rightarrow K$ is a compact homotopy and $0\not\in(I-H)(\partial\,U\times[0,1])$, then
		\[i_{K}(H(\cdot,0),U)=i_{K}(H(\cdot,1),U).\]
		\item (Normalization) If $N$ is a constant map with $Nx=\bar{x}$ for every $x\in\overline{U}$, then
		\[i_{K}(N,U)=\left\{\begin{array}{ll} 1, & \text{ if } \bar{x}\in U, \\ 0, & \text{ if } \bar{x}\not\in\overline{U}. \end{array} \right. \]
	\end{enumerate}
\end{proposition}

The following computations of the fixed point index will be also useful in our reasoning. Their proofs can be found, for instance, in \cite[Lemma 2.3.1 and 2.3.2]{guolak}.

\begin{proposition}\label{prop_ind01}
	Let $K$ be a cone, $U\subset K$ be a bounded relatively open set such that $0\in U$ and $N:\overline{U}\rightarrow K$ be a compact map without fixed points on $\partial_{K}\,U$.
	\begin{enumerate}[$(a)$]
		\item If $Nx\neq \lambda\, x$ for all $x\in \partial_{K}\, U$ and all $\lambda> 1$, then $i_{K}(N,U)=1$.
		\item If there exists a compact map $S:\overline{U}\rightarrow K$ such that 
		\[\inf_{x\in \overline{U}}\left\|S x\right\|>0, \text{ and} \]
		\[x-N x\neq \mu\, S x \ \text{ for all } x\in \partial_{K}\, U \text{ and every } \mu> 0 ,\] 
		then $i_{K}(N,U)=0$.
	\end{enumerate}	
\end{proposition}

In case of operators defined in the Cartesian product of normed spaces, a component-wise combination of the previous conditions is possible. The computation of the fixed point index in such situation, which is detailed in our next proposition, plays a key role in the proof of the main results of the following section. 

\begin{proposition}\label{prop_ind_sys}
	Let $U\times V$ be a bounded relatively open subset of the cone $K= K_1\times K_2$ in the normed spaces product $X=X_1\times X_2$, such that $0\in U$. Assume that $N=(N_1,N_2):\overline{U\times V}\rightarrow K$ and $S:\overline{U\times V}\rightarrow K_2$ are compact mappings satisfying the following conditions:
	\begin{enumerate}[$(a)$]
		\item $N_1 x \neq \lambda\, x_1$ for all $x_1\in \partial_{K_1}U$, $x_2\in\overline{V}$ and all $\lambda>1$;
		\item \begin{enumerate}[$(i)$]
			\item $\inf_{x\in \overline{U\times V}}\left\|S x\right\|>0$;
			\item $x_2-N_2 x\neq \mu\, S x$ for all $x_1\in\overline{U}$, $x_2\in\partial_{K_2} V$ and every $\mu>0$.
		\end{enumerate}
	\end{enumerate}	
	If $N$ has no fixed points on $\partial_K\, (U\times V)$, then $i_{K}(N,U\times V)=0$.
\end{proposition}

\noindent
{\bf Proof.} First, let us define the following positive numbers
\[\alpha:=\sup_{x_2\in\overline{V}}\left\|x_2\right\|, \quad \beta:=\sup_{x\in\overline{U\times V}}\left\|N_2 x\right\| \quad \text{and} \quad \gamma:=\inf_{x\in \overline{U\times V}}\left\|Sx\right\|. \]
The existence of $\alpha$ and $\beta$ is guaranteed by the fact that $\overline{U}$ and $\overline{V}$ are bounded sets and $N$ is a compact map. Moreover, by hypothesis $(b)$-$(i)$, we have $\gamma>0$. 

Now, let us fix a positive number $\mu_0>(\alpha+\beta)/\gamma$ and consider the compact homotopy $H:\overline{U\times V}\times [0,1]\rightarrow K$ given by 
\[H(x,t)=\left(t\,N_1 x,N_2 x+(1-t)\mu_0\,Sx\right). \]
It is an \textit{admissible} homotopy, that is, $H$ has no fixed points on $\partial_K\,(U\times V)\times [0,1]$. Otherwise, there exist $(\bar{x}_1,\bar{x}_2)\in \left(\partial_{K_1} U\times \overline{V}\right)\cup \left(\overline{U}\times \partial_{K_2} V\right)$ and $\bar{t}\in [0,1]$ such that
\[\bar{x}_1=\bar{t}\,N_1(\bar{x}_1,\bar{x}_2), \quad \bar{x}_2=N_2(\bar{x}_1,\bar{x}_2)+(1-\bar{t})\mu_0\,S(\bar{x}_1,\bar{x}_2). \]
Since $N$ has no fixed points on $\partial_K\, (U\times V)$, it follows that $\bar{t}<1$. If $\bar{t}=0$, then $\bar{x}_1=0\in U$, $\bar{x}_2\in\partial_{K_2}V$ and $\bar{x}_2=N_2(0,\bar{x}_2)+(1-\bar{t})\mu_0\,S(0,\bar{x}_2),$ a contradiction with $(b)$-$(ii)$. Hence, $\bar{t}\in (0,1)$ and so either \[\dfrac{1}{\bar{t}}\,\bar{x}_1=N_1(\bar{x}_1,\bar{x}_2) \] contradicts $(a)$ if $(\bar{x}_1,\bar{x}_2)\in \partial_{K_1} U\times \overline{V}$ or $\bar{x}_2-N_2(\bar{x}_1,\bar{x}_2)=(1-\bar{t})\mu_0\,S(\bar{x}_1,\bar{x}_2)$ contradicts $(b)$-$(ii)$ in case that $(\bar{x}_1,\bar{x}_2)\in \overline{U}\times \partial_{K_2} V$. 

Therefore, by the homotopy invariance of the fixed point index, we obtain
\[i_{K}(N,U\times V)=i_{K}(H(\cdot,0),U\times V). \]
Finally, observe that the map $H(\cdot,0)$ has no fixed points in $U\times V$. Indeed, if there exists such a fixed point $(\bar{x}_1,\bar{x}_2)$, then
\[\bar{x}_1=0, \quad \bar{x}_2=N_2(0,\bar{x}_2)+\mu_0\,S(0,\bar{x}_2). \]
From the second equality, it follows that
\[\mu_0=\dfrac{\left\|\bar{x}_2-N_2(0,\bar{x}_2)\right\|}{\left\|S(0,\bar{x}_2)\right\|}\leq \dfrac{\alpha+\beta}{\gamma}, \]
a contradiction with the choice of $\mu_0$. In conclusion, by the existence property of the index, it is deduced that $i_{K}(H(\cdot,0),U\times V)=0$. 
\qed

\begin{remark}
	It is clear that the conclusion of Proposition \ref{prop_ind_sys} remains valid if the roles of $N_1$ and $N_2$ expressed by assumptions $(a)$ and $(b)$ are interchanged. 
\end{remark}

\begin{remark}
	Note that Proposition \ref{prop_ind_sys} is a refinement of \cite[Lemma 2.3]{PreRod} even if $S$ is a constant map such that $Sx=h\in K_{2}\setminus\{0\}$ for every $x\in\overline{U\times V}$. Indeed, the corresponding result from \cite{PreRod} requires that condition $(a)$ (respectively $(b)$ with $S\equiv h$) also hold for $\lambda= 1$ (resp. for $\mu=0$), which is stronger than assuming that $N$ is fixed point free on the boundary of $U\times V$.
\end{remark}

\section{A fixed point theorem for operator systems}\label{sec_main}
Let $(X_1,\left\|\cdot\right\|_{X_1})$ and $(X_2,\left\|\cdot\right\|_{X_2})$ be normed linear spaces and $K_1\subset X_1$ and $K_2\subset X_2$ two cones. We will employ the notation $K:=K_1\times K_2$ for the corresponding cone in the cartesian product $X:=X_1\times X_2$. Notice that $X$ is also a normed space endowed with the norm $\left\|(x_1,x_2)\right\|_{X}:=\max\{\left\|x_1\right\|,\left\|x_2\right\| \}$. When no confusion may occur, all the norms $\left\|\cdot\right\|_{X_1}$, $\left\|\cdot\right\|_{X_2}$ and $\left\|\cdot\right\|_{X}$ will be simply denoted as $\left\|\cdot\right\|$.

Given $\tau_1,\tau_2>0$ and $\tau=(\tau_1,\tau_2)$, we will use the following notations:
\[K_{\tau}:=\left\{x=(x_1,x_2)\in K: \left\|x_i\right\|< \tau_i \text{ for } i=1,2 \right\}=(K_1)_{\tau_1}\times (K_2)_{\tau_2}, \text{ and} \]
\[\overline{K}_{\tau}:=\left\{x=(x_1,x_2)\in K: \left\|x_i\right\|\leq \tau_i \text{ for } i=1,2 \right\}=(\overline{K}_1)_{\tau_1}\times (\overline{K}_2)_{\tau_2},\]
where for each $i\in\{1,2\}$,
\[(K_i)_{\tau_i}:=\left\{x_i\in K_i:\left\|x_i\right\|< \tau_i \right\} \]
and $(\overline{K}_i)_{\tau_i}$ denotes its closure. Now, for $0<r_i<R_i$ ($i=1,2$) fixed, we consider the following sets:
\[K_{r,R}:=\left\{x=(x_1,x_2)\in K: r_i<\left\|x_i\right\|< R_i \text{ for } i=1,2 \right\}=(K_1)_{r_1,R_1}\times (K_2)_{r_2,R_2}, \text{ and} \] 
\[\overline{K}_{r,R}:=\left\{x=(x_1,x_2)\in K: r_i\leq\left\|x_i\right\|\leq R_i \text{ for } i=1,2\right\}=(\overline{K}_1)_{r_1,R_1}\times (\overline{K}_2)_{r_2,R_2}, \]
with $(K_i)_{r_i,R_i}:=(K_i)_{R_i}\setminus(\overline{K}_i)_{r_i}$ and $(\overline{K}_i)_{r_i,R_i}:=(\overline{K}_i)_{R_i}\setminus({K}_i)_{r_i}$, $i=1,2$.

\begin{figure}[h]
	\centering
	\includegraphics[scale=0.95]{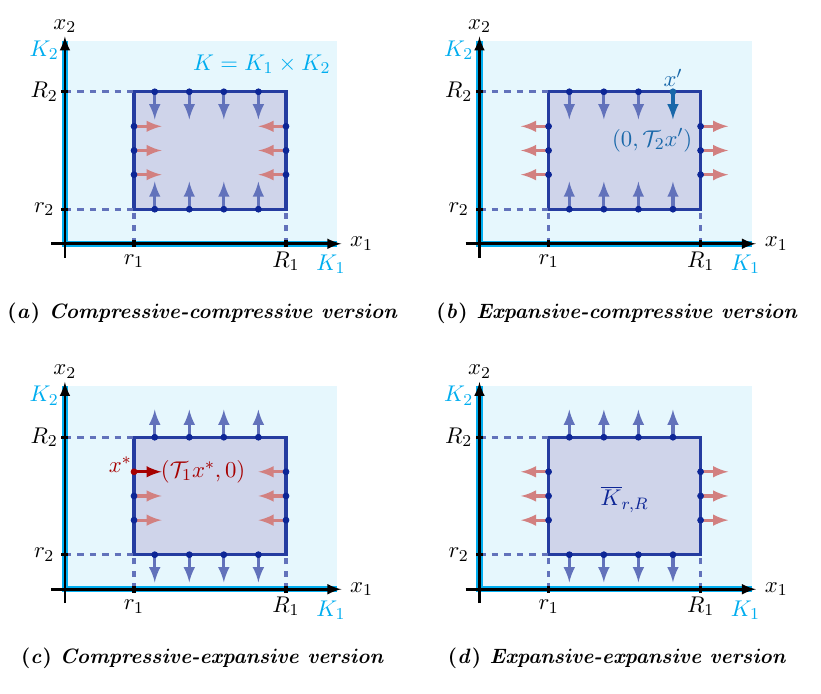}
	\caption{Visualization of all possible cases in Theorem \ref{th}. }
	\label{fig}
\end{figure}

The primary goal of this section is to prove our main result, Theorem \ref{main_th}. It allows the operator to exhibit different behavior in each of its components (either compression or expansion), as described below. See also Figure~\ref{fig}.

\begin{theorem}
	\label{th}
	Let $r,R\in\mathbb{R}^2$ with $0<r_i<R_i$ ($i=1,2$) and assume that $\mathcal{T}=(\mathcal{T}_1,\mathcal{T}_2):\overline{K}_{r,R}\rightarrow K$ is a compact map such that one of the following conditions is satisfied in $\overline{K}_{r,R}$:
	\begin{enumerate}[$(a)$]
		\item \textbf{Compressive-compressive version}
		\begin{enumerate}[$(i)$]
			\item $\|\mathcal{T}_1 x\|\geq\|x_1\|$ if $\|x_1\|=r_1$ and $\|\mathcal{T}_1 x\|\leq\|x_1\|$ if  $\|x_1\|=R_1$;
			\item $\|\mathcal{T}_2 x\|\geq\|x_2\|$ if $\|x_2\|=r_2$ and $\|\mathcal{T}_2 x\|\leq\|x_2\|$ if  $\|x_2\|=R_2$.
		\end{enumerate}
		\item \textbf{Expansive-compressive version}
		\begin{enumerate}[$(i)$]
			\item $\|\mathcal{T}_1 x\|\leq\|x_1\|$ if $\|x_1\|=r_1$ and $\|\mathcal{T}_1 x\|\geq\|x_1\|$ if  $\|x_1\|=R_1$;
			\item $\|\mathcal{T}_2 x\|\geq\|x_2\|$ if $\|x_2\|=r_2$ and $\|\mathcal{T}_2 x\|\leq\|x_2\|$ if  $\|x_2\|=R_2$.
		\end{enumerate}
		
		\item \textit{\textbf{Compressive-expansive version}}
		\begin{enumerate}[$(i)$]
			\item $\|\mathcal{T}_1 x\|\geq\|x_1\|$ if $\|x_1\|=r_1$ and $\|\mathcal{T}_1 x\|\leq\|x_1\|$ if $\|x_1\|=R_1$;
			\item $\|\mathcal{T}_2 x\|\leq\|x_2\|$ if $\|x_2\|=r_2$ and $\|\mathcal{T}_2 x\|\geq\|x_2\|$ if $\|x_2\|=R_2$.
		\end{enumerate}
		\item \textit{\textbf{Expansive-expansive version}}
		\begin{enumerate}[$(i)$]
			\item $\|\mathcal{T}_1 x\|\leq\|x_1\|$ if $\|x_1\|=r_1$ and $\|\mathcal{T}_1 x\|\geq\|x_1\|$ if $\|x_1\|=R_1$;
			\item $\|\mathcal{T}_2 x\|\leq\|x_2\|$ if $\|x_2\|=r_2$ and $\|\mathcal{T}_2 x\|\geq\|x_2\|$ if $\|x_2\|=R_2$.
		\end{enumerate}
	\end{enumerate}
	Then $\mathcal{T}$ has at least one fixed point $x=(x_1,x_2)\in K$ such that $r_i\leq \left\|x_i\right\|\leq R_i$ ($i=1,2$).	
\end{theorem}

To demonstrate the particular case of Theorem \ref{th} in which both components of the operator are of compression type, we compute the fixed point index of an extension of the map $\mathcal{T}$ over different open sets. Let us define this compact extension of $\mathcal{T}$. First, we observe that $\overline{K}_{r,R}$ is a retract of $\overline{K}_{R}$. Indeed, following \cite[Example 3]{fel} or \cite{JRL}, we can define the retraction $\rho=(\rho_1,\rho_2):\overline{K}_{R}\rightarrow \overline{K}_{r,R}$ as
\begin{equation}\label{ret}
	\rho_i(x_i)=\left\{\begin{array}{ll} r_i\dfrac{x_i+(r_i-\left\|x_i\right\|)^2 h_i}{\left\|x_i+(r_i-\left\|x_i\right\|)^2 h_i\right\|} & \text{ if } \left\|x_i\right\|<r_i, \\ x_i & \text{ if } r_i\leq \left\|x_i\right\|\leq R_i, \end{array} \right.
\end{equation}
where $h_i\in K_i\setminus\{0\}$ is fixed, $i=1,2$.

This retraction $\rho$ allows us to extend the definition of a map $\mathcal{T}$ under the assumptions of Theorem~\ref{th}~(a) to the set $\overline{K}_{R}$ simply as:
\begin{equation}\label{eq_N}
	N=(N_1,N_2):\overline{K}_{R}\rightarrow K, \quad N:=\mathcal{T}\circ\rho.
\end{equation}
Since $\mathcal{T}$ is a compact map, then so is $N$. In addition, the following compression type behavior for both components of $N$ is fulfilled in the set $\overline{K}_{R}$:
\begin{equation}\label{eq_N_comp}
	\left\|N_i x\right\|\geq \left\|x_i\right\| \ \text{ if } \left\|x_i\right\|=r_i \text{ and } \left\|N_i x\right\|\leq \left\|x_i\right\| \text{ if } \left\|x_i\right\|=R_i \quad (i=1,2).
\end{equation}
Let us check it for $i=1$ (the case $i=2$ is analogous). Take $x=(x_1,x_2)\in \overline{K}_{R}$ such that $\left\|x_1\right\|=r_1$. By definition of $N$ and $\rho$, we have $N_1(x_1,x_2)=\mathcal{T}_1(\rho_1(x_1),\rho_2(x_2))=\mathcal{T}_1(x_1,\rho_2(x_2))$. Observe that $(x_1,\rho_2(x_2))\in \overline{K}_{r,R}$ and $\left\|x_1\right\|=r_1$, which thanks to Theorem \ref{th} $(a)$-$(i)$ implies $\left\|N_1(x_1,x_2)\right\|=\left\|\mathcal{T}_1(x_1,\rho_2(x_2))\right\|\geq \left\|x_1\right\|.$
Now, choose $x=(x_1,x_2)\in \overline{K}_{R}$ with $\left\|x_1\right\|=R_1$. It follows, due to the second part of the same condition, that $\left\|N_1(x_1,x_2)\right\|=\left\|\mathcal{T}_1(x_1,\rho_2(x_2))\right\|\leq \left\|x_1\right\|,$ and, thus, \eqref{eq_N_comp} is proven.

Next result ensures that the fixed point index of the compact map $N$ with respect to the cone $K$ is well-defined over the open sets $K_{r}$, $K_{R}$, $K_{r,R}$, $(K_1)_{R_1}\times (K_2)_{r_2}$ and $(K_1)_{r_1}\times (K_2)_{R_2}$.

\begin{lemma}\label{lem_NnoFP}
	Let $\mathcal{T}$ be a map under the hypotheses of Theorem \ref{th} (\textit{a}). 
	
	If $\mathcal{T}$ is fixed point free on the boundary of the set $\overline{K}_{r,R}$, then its extension $N$ defined as \eqref{eq_N} also has no fixed point $x=(x_1,x_2)\in \overline{K}_R$ such that $\left\|x_i\right\|=r_i$ or $\left\|x_i\right\|=R_i$ with $i\in\{1,2\}$.	
\end{lemma}

\noindent
{\bf Proof.} Let us suppose that $N$ has a fixed point $(x_1,x_2)\in \overline{K}_{R}$ such that $\left\|x_1\right\|\in \{r_1,R_1\}$. Since $N=\mathcal{T}$ in $\overline{K}_{r,R}$ and $\mathcal{T}$ has no fixed points on the boundary of $\overline{K}_{r,R}$, it follows that $\left\|x_2\right\|<r_2$. Then
\[(x_1,x_2)=N(x_1,x_2)=\left(\mathcal{T}_1(x_1,\rho_2(x_2)),\mathcal{T}_2(x_1,\rho_2(x_2))\right).\]
 
By \eqref{ret} and the fact that $\left\|x_2\right\|<r_2$, it is clear that $\left\|\rho_2(x_2)\right\|=r_2$. Hence, the cone-compression condition given by Theorem \ref{th} $(a)$-$(ii)$ implies that 
\[\left\|x_2\right\|=\left\|\mathcal{T}_2(x_1,\rho_2(x_2))\right\|\geq \left\|\rho_2(x_2)\right\|=r_2, \]
a contradiction.

Therefore, it is proven that $N$ has no a fixed point $(x_1,x_2)\in \overline{K}_{R}$ with $\left\|x_1\right\|=r_1$ or $\left\|x_1\right\|=R_1$. Similarly, it can be shown that $N$ does not have a fixed point $(x_1,x_2)$ such that $\left\|x_2\right\|=r_2$ or $\left\|x_2\right\|=R_2$ by interchanging the role of the maps $\mathcal{T}_1$ and $\mathcal{T}_2$ in the previous reasoning.
\qed

Now, we are in a position to establish the following computation of the fixed point index of the map $\mathcal{T}$ in \textit{Compressive-compressive case} of Theorem \ref{th}.

\begin{theorem}\label{th_index_comp}
	Let $r,R\in\mathbb{R}^2$ with $0<r_i<R_i$ ($i=1,2$) and assume that $\mathcal{T}=(\mathcal{T}_1,\mathcal{T}_2):\overline{K}_{r,R}\rightarrow K$ is a compact map such that for each $i\in\{1,2\}$ the following condition is satisfied in $\overline{K}_{r,R}$:
		\begin{equation}\label{eq_comp_T}
		\left\|\mathcal{T}_i x\right\|\geq \left\|x_i\right\| \ \text{ if } \left\|x_i\right\|=r_i \text{ and } \left\|\mathcal{T}_i x\right\|\leq \left\|x_i\right\| \text{ if } \left\|x_i\right\|=R_i.
	\end{equation}
	If $\mathcal{T}$ has no fixed points on the boundary of $\overline{K}_{r,R}$, then $i_{K}(\mathcal{T},K_{r,R})=1$.
\end{theorem}

\noindent
{\bf Proof.} Firstly, consider the map $N$ as defined in \eqref{eq_N}. Let us compute the fixed point index of $N$ on different open sets, namely, $K_{R}$, $K_{r}$, $(K_1)_{R_1}\times (K_2)_{r_2}$, $(K_1)_{r_1}\times (K_2)_{R_2}$ and $(K_1)_{r_1,R_1}\times (K_2)_{r_2}$. Observe that Lemma \ref{lem_NnoFP} ensures the index of $N$ on these sets is well-defined. Moreover, condition \eqref{eq_comp_T} guarantees that $N$ satisfies \eqref{eq_N_comp}, as already shown.

Note that condition 
\[\left\|N_i x\right\|\leq \left\|x_i\right\| \text{ if } \left\|x_i\right\|=R_i \quad (i=1,2)\]
ensures that
\begin{equation}\label{eq_Ni_Ri}
	N_i x\neq \lambda\, x_i \quad \text{ for all } x\in \overline{K}_{R} \text{ with } \left\|x_i\right\|=R_i \text{ and all } \lambda>1 \quad (i=1,2).
\end{equation}
If not, there exist $(\bar{x}_1,\bar{x}_2)\in \overline{K}_{R}$  with $\left\|\bar{x}_i\right\|=R_i$ and $\bar{\lambda}>1$ such that $N_i \bar{x}=\bar{\lambda}\,\bar{x}_i$. Thus,
\[\left\|N_i \bar{x}\right\|=\left\|\bar{\lambda}\,\bar{x}_i\right\|=\bar{\lambda}\left\|\bar{x}_i\right\|>\left\|\bar{x}_i\right\|, \]
a contradiction. By \eqref{eq_Ni_Ri}, it follows that $N x\neq \lambda\,x$ for all $x\in \partial_{K} K_{R}$ and all $\lambda>1$. Then Proposition~\ref{prop_ind01}~$(a)$ yields 
\[i_{K}(N,K_R)=1. \]

Similarly, it is easy to check that condition 
\[\left\|N_i x\right\|\geq \left\|x_i\right\| \text{ if } \left\|x_i\right\|=r_i \quad (i=1,2)\]
implies that
\begin{equation}\label{eq_Ni_ri}
	N_i x\neq \mu\, x_i \quad \text{ for all } x\in \overline{K}_{R} \text{ with } \left\|x_i\right\|=r_i \text{ and all } \mu\in (0,1) \quad (i=1,2).
\end{equation}
Condition \eqref{eq_Ni_ri} implies, in turn, the following one:
\begin{equation}\label{eq_Ni_ri2}
	x_i-N_i x\neq \mu\,N_i x \quad \text{ for all } x\in \overline{K}_{R} \text{ with } \left\|x_i\right\|=r_i \text{ and all } \mu>0 \quad (i=1,2).  
\end{equation}
Indeed, if \eqref{eq_Ni_ri2} does not hold, there exist $(\bar{x}_1,\bar{x}_2)\in \overline{K}_{R}$  with $\left\|\bar{x}_i\right\|=r_i$ and $\bar{\mu}>0$ such that $\bar{x}_i-N_i \bar{x}=\bar{\mu}\,N_i \bar{x}$ and so $N_i \bar{x}=(1+\bar{\mu})^{-1}\bar{x}_i$ with $(1+\bar{\mu})^{-1}\in (0,1)$, which contradicts \eqref{eq_Ni_ri}.

Due to \eqref{eq_Ni_ri2} and defining $S:=N$ in $\overline{K}_r$, we have that 
\[x-N x\neq \mu\, S x \ \text{ for all } x\in \partial_{K}\, K_r \text{ and every } \mu> 0. \]
Moreover,
\[\inf_{x\in \overline{K}_r}\left\|Sx\right\|=\inf_{x\in \overline{K}_r}\left\|\mathcal{T}(\rho(x))\right\|_{X_1\times X_2}\geq \inf_{x\in \overline{K}_r}\left\|\mathcal{T}_1(\rho_1(x_1),\rho_2(x_2))\right\|_{X_1}\geq \left\|\rho_1(x_1)\right\|=r_1, \]
where the last inequality follows from \eqref{eq_comp_T}, since the definition of the retraction $\rho$ ensures $\left\|\rho_1(x_1)\right\|=r_1$. 
Therefore, Proposition~\ref{prop_ind01}~$(b)$ yields
\[i_{K}(N,K_r)=0. \]

Now, let us study the index of $N$ on the open set $(K_1)_{R_1}\times (K_2)_{r_2}$. By \eqref{eq_Ni_Ri} and \eqref{eq_Ni_ri2}, we have
\begin{align*}
		N_1 x\neq \lambda\, x_1 \quad \text{ for all } x_1\in \partial_{K_1}(K_1)_{R_1}, \ x_2\in(\overline{K}_2)_{r_2} \text{ and all } \lambda>1, \\	
		x_2-N_2 x\neq \mu\,\tilde{S} x \quad \text{ for all } x_1\in (\overline{K}_{1})_{R_1}, \  x_2\in\partial_{K_2}(K_2)_{r_2} \text{ and all } \mu>0,
\end{align*}
with $\tilde{S}:=N_2$. Moreover,
\[\inf_{x\in (\overline{K}_1)_{R_1}\times (\overline{K}_2)_{r_2}}\|\tilde{S} x\|=\inf_{x\in (\overline{K}_1)_{R_1}\times (\overline{K}_2)_{r_2}}\left\|\mathcal{T}_2(\rho_1(x_1),\rho_2(x_2))\right\|\geq \left\|\rho_2(x_2)\right\|=r_2>0.  \]
Hence, Proposition \ref{prop_ind_sys} gives
\[i_{K}(N, (K_1)_{R_1}\times (K_2)_{r_2})=0. \]
Similarly, it can be seen that $i_{K}(N, (K_1)_{r_1}\times (K_2)_{R_2})=0$.

Finally, the conclusion is a consequence of the additivity property of the fixed point index since
\[i_{K}(N,(K_1)_{r_1,R_1}\times (K_2)_{r_2})=i_{K}(N, (K_1)_{R_1}\times (K_2)_{r_2})-i_{K}(N,K_r)=0 \] 
and thus
\[i_{K}(N,K_{r,R})=i_{K}(N,K_R)-i_{K}(N,(K_1)_{r_1,R_1}\times (K_2)_{r_2})-i_{K}(N, (K_1)_{r_1}\times (K_2)_{R_2})=1.\]
By definition, $N=\mathcal{T}$ in $\overline{K}_{r,R}$, so $i_{K}(\mathcal{T},K_{r,R})=1$.
\qed

\begin{remark}
	Clearly, Theorem \ref{th_index_comp} implies the existence of a fixed point of $\mathcal{T}$, as stated in Theorem \ref{th}. Indeed, either $\mathcal{T}$ has a fixed point on the boundary of $\overline{K}_{r,R}$ or $i_{K}(\mathcal{T},K_{r,R})=1\neq 0$, which ensures the existence of a fixed point of $\mathcal{T}$ in $K_{r,R}$, thanks to the existence property of the fixed point index.
\end{remark}


The proof of Theorem \ref{th} for the rest of situations is based on a trick to transform cone-expansive into cone-compressive conditions, exactly as in \cite{PrecupFPT,PrecupSDC}. Some details are included for the sake of completeness.

\noindent
{\bf Proof of Theorem \ref{th}.}
In \textit{Expansive-compressive case} if we consider the map $\mathcal{T}^*_1:\overline{K}_{r,R}\longrightarrow K_1$ given by 
\begin{equation*}
	\mathcal{T}^*_1 x=\left(\frac{R_1}{\|x_1\|}+\frac{r_1}{\|x_1\|}-1\right)^{-1}\mathcal{T}_1\left(\left(\frac{R_1}{\|x_1\|}+\frac{r_1}{\|x_1\|}-1\right)x_1,x_2\right).
\end{equation*}
%
%
It is easy to check that $\tilde{\mathcal{T}}=(\mathcal{T}^*_1,\mathcal{T}_2)$ is under the conditions of \textit{(a)} and thus $\tilde{\mathcal{T}}$ has a fixed point in $\overline{K}_{r,R}$. Let $(v_1,v_2)$ be such a fixed point, then the point $(u_1,u_2)\in\overline{K}_{r,R}$ given by 
\begin{equation*}
	u_1=\left(\frac{R_1}{\|x_1\|}+\frac{r_1}{\|x_1\|}-1\right)v_1 \quad \text{ and }\quad u_2=v_2,
\end{equation*}
is a fixed point of $\mathcal{T}$.

The \textit{Compressive-expansive case} is analogous, now with the operator $\hat{\mathcal{T}}=(\mathcal{T}_1,\mathcal{T}^*_2)$, where $\mathcal{T}^*_2:\overline{K}_{r,R}\longrightarrow K_2$ is defined as
\begin{equation*}
	\mathcal{T}^*_2 x=\left(\frac{R_2}{\|x_2\|}+\frac{r_2}{\|x_2\|}-1\right)^{-1}\mathcal{T}_2\left(x_1,\left(\frac{R_2}{\|x_2\|}+\frac{r_2}{\|x_2\|}-1\right)x_2\right).
\end{equation*}
%

Finally, for \textit{Expansive-expansive case}, we consider the auxiliary operator $\mathcal{T}^*=(\mathcal{T}^*_1,\mathcal{T}^*_2)$, which satisfies the compressive-compressive version and so it has a fixed point. If $(v_1,v_2)\in\overline{K}_{r,R}$ is a fixed point of $\mathcal{T}^*$, then the point $(u_1,u_2)\in\overline{K}_{r,R}$ given by
\begin{equation*}
	u_1=\left(\frac{R_1}{\|x_1\|}+\frac{r_1}{\|x_1\|}-1\right)v_1 \quad \text{ and }\quad u_2=\left(\frac{R_2}{\|x_2\|}+\frac{r_2}{\|x_2\|}-1\right)v_2,
\end{equation*}
is a fixed point of $\mathcal{T}$.
\qed



Observe that if both components of the operator $\mathcal{T}=(\mathcal{T}_1,\mathcal{T}_2)$ are of compression type, then it is not only shown that $\mathcal{T}$ has a fixed point, but its fixed point index in $K_{r,R}$ is equal to $1$. This additional information does have interesting consequences, for instance in order to obtain multiple fixed points of the operator.

Nevertheless, the value of the index is not calculated if some of the components is expansive. The reason is that, in such a situation, it is not possible to ensure that the extension $N$ has no fixed points on the boundaries of the relevant open sets, solely assuming that $\mathcal{T}$ does not have them on $\partial_{K} K_{r,R}$, as in Lemma \ref{lem_NnoFP}. 

If we strengthen the hypotheses of Theorem \ref{main_th} by replacing the non-strict inequalities with strict inequalities, the value of the fixed point index is obtained in each case.   

\begin{theorem}
	\label{th_cond_without_ec}
	Let $r,R\in\mathbb{R}^2$ with $0<r_i<R_i$ ($i=1,2$) and assume that $\mathcal{T}=(\mathcal{T}_1,\mathcal{T}_2):\overline{K}_{r,R}\rightarrow K$ is a compact map and for each $i\in\{1,2\}$ either of the following conditions holds in $\overline{K}_{r,R}$:
	\begin{enumerate}
		\item[(a)] $\left\|\mathcal{T}_i x\right\|> \left\|x_i\right\|$ if $\left\|x_i\right\|=r_i$ and  $\left\|\mathcal{T}_i x\right\|< \left\|x_i\right\|$ if $\left\|x_i\right\|=R_i$; or
		\item[(b)] $\left\|\mathcal{T}_i x\right\|< \left\|x_i\right\|$ if $\left\|x_i\right\|=r_i$ and  $\left\|\mathcal{T}_i x\right\|> \left\|x_i\right\|$ if $\left\|x_i\right\|=R_i$. 
	\end{enumerate}
	Then \[i_{K}(\mathcal{T},K_{r,R})=(-1)^s, \]
	where $s\in\{0,1,2\}$ denotes the number of expansive components of $\mathcal{T}$. 
\end{theorem}

\noindent
{\bf Proof.} Suppose that $\mathcal{T}_1$ satisfies condition $(a)$ and $\mathcal{T}_2$, condition $(b)$ (the remaining cases are analogous). 
Consider the operator $N$ defined as in \eqref{eq_N}. It fulfills the following inequalities in $\overline{K}_{R}$:
\begin{align*}
	&\left\|N_1 x\right\|> \left\|x_1\right\| \text{ if } \left\|x_1\right\|=r_1, \quad \left\|N_1 x\right\|< \left\|x_1\right\| \text{ if } \left\|x_1\right\|=R_1; \text{ and } \\
	&\left\|N_2 x\right\|< \left\|x_2\right\| \text{ if } \left\|x_2\right\|=r_2, \quad \left\|N_2 x\right\|> \left\|x_2\right\| \text{ if } \left\|x_2\right\|=R_2.
\end{align*} 
Notice that the strict inequalities clearly ensure that $N$ has no fixed points on the boundaries of the open sets $K_{R}$, $K_{r}$, $(K_1)_{R_1}\times (K_2)_{r_2}$ and $(K_1)_{r_1}\times (K_2)_{R_2}$. Therefore, following the reasoning in the proof of Theorem \ref{th_index_comp}, it is deduced from Proposition \ref{prop_ind01} that $i_{K}(N,(K_1)_{R_1}\times (K_2)_{r_2} )=1$ and $i_{K}(N, (K_1)_{r_1}\times (K_2)_{R_2})=0.$ Moreover, by Proposition \ref{prop_ind_sys}, we also have $i_{K}(N,K_r)=0=i_{K}(N,K_R).$
	
Finally, applying the additivity property of the fixed point index twice as in the mentioned proof, we obtain
\[i_{K}(N,(K_1)_{r_1,R_1}\times (K_2)_{r_2} )=i_{K}(N,(K_1)_{R_1}\times (K_2)_{r_2} )-i_{K}(N,K_r)=1 \]
and then
\[i_{K}(N,K_{r,R})=i_{K}(N,K_{R})-i_{K}(N,(K_1)_{r_1,R_1}\times (K_2)_{r_2} )-i_{K}(N, (K_1)_{r_1}\times (K_2)_{R_2})=-1. \]
Therefore, $i_{K}(\mathcal{T},K_{r,R})=-1$, as wished.
\qed

Next multiplicity result, in the line of Amann \cite{amann}, can be deduced as a straightforward consequence of the previous fixed point index computation.

\begin{theorem}
Let $r^{(j)},R^{(j)}\in\mathbb{R}^2$ with $0<r_i^{(j)}<R_i^{(j)}$ $(i=1,2,\ j=1,2,\dots,2m+1, \ m\in\mathbb{N})$. Assume that the sets $\overline{K}_{r^{(j)},R^{(j)}}$ are such that
\[\bigcup_{j=1}^{2m}\overline{K}_{r^{(j)},R^{(j)}}\subset\overline{K}_{r^{(2m+1)},R^{(2m+1)}}\quad \text{and}\quad \overline{K}_{r^{(k)},R^{(k)}}\cap \overline{K}_{r^{(l)},R^{(l)}}=\emptyset \quad k,l\in\{1,\dots,2m\}, \ k\neq l.\]

Moreover, assume that $\mathcal{T}=(\mathcal{T}_1,\mathcal{T}_2):\overline{K}_{r^{(2m+1)},R^{(2m+1)}}\rightarrow K$ is a compact map and for each $i\in\{1,2\}$ and each $j\in\{1,\dots,2m+1\}$ one of the following conditions is satisfied in $\overline{K}_{r^{(j)},R^{(j)}}$:
\begin{enumerate}
	\item[(a)] $\left\|\mathcal{T}_i x\right\|> \left\|x_i\right\|$ if $\left\|x_i\right\|=r_i^{(j)}$ and  $\left\|\mathcal{T}_i x\right\|< \left\|x_i\right\|$ if $\left\|x_i\right\|=R_i^{(j)}$; or
	\item[(b)] $\left\|\mathcal{T}_i x\right\|< \left\|x_i\right\|$ if $\left\|x_i\right\|=r_i^{(j)}$ and  $\left\|\mathcal{T}_i x\right\|> \left\|x_i\right\|$ if $\left\|x_i\right\|=R_i^{(j)}$. 
\end{enumerate}
Then $\mathcal{T}$ has at least $2m+1$ distinct fixed points $\bar{x}^1, \bar{x}^2, \dots, \bar{x}^{2m+1}$ such that
\[\bar{x}^j\in K_{r^{(j)},R^{(j)}} \quad (j=1,2\dots,2m)\, \text{ and }\, \bar{x}^{2m+1}\in K_{r^{(2m+1)},R^{(2m+1)}}\backslash\bigcup_{j=1}^{2m}\overline{K}_{r^{(j)},R^{(j)}}.\] 
\end{theorem}

\noindent
{\bf Proof.} It follows from Theorem \ref{th_cond_without_ec} that
\begin{equation}
	\label{ec3.8}
i_K(\mathcal{T},K_{r^{(j)},R^{(j)}})=\pm 1 \quad (j=1,2,\dots,2m+1).
\end{equation}
Since it does not vanish, one obtains that $\mathcal{T}$ has fixed points $\bar{x}^j\in K_{r^{(j)},R^{(j)}}$ with $j=1,\dots,2m$. These fixed points are different because the sets $K_{r^{(j)},R^{(j)}}$ are pairwise disjoint.  

Next, the additivity property of the fixed point index, together with \eqref{ec3.8}, implies
\[i_K\left(\mathcal{T},K_{r^{(2m+1)},R^{(2m+1)}}\backslash\bigcup_{j=1}^{2m}\overline{K}_{r^{(j)},R^{(j)}}\right)=i_K(\mathcal{T},K_{r^{(2m+1)},R^{(2m+1)}})-\sum_{j=1}^{2m}i_K(\mathcal{T},K_{r^{(j)},R^{(j)}})\]
is an odd number and so it is different from zero.
Therefore, the operator $\mathcal{T}$ has an additional fixed point $\bar{x}^{2m+1}\in K_{r^{(2m+1)},R^{(2m+1)}}\backslash\bigcup_{j=1}^{2m}\overline{K}_{r^{(j)},R^{(j)}}$. \qed

\section{Applications}\label{sec_app}

In this section, we deal with the existence of positive solutions for the following system of second-order equations
\begin{equation}\label{eq_periodic_sys}
	\begin{array}{ll} x''+a_1(t)\,x=f_1(t,x,y) & \quad t\in[0,T], \\ y''+a_2(t)\,y=f_2(t,x,y) & \quad t\in[0,T], \\ x(0)-x(T)=0=x'(0)-x'(T), & \\ y(0)-y(T)=0=y'(0)-y'(T), \end{array} 
\end{equation}
where $a_i\in L^{p_i}(0,T)$ with $1\leq p_i\leq +\infty$, $i=1,2$, and $f_1,f_2:[0,T]\times \mathbb{R}_{+}^2\rightarrow\mathbb{R}$ are continuous functions.

First, we assume that for each $i\in\{1,2\}$ the following condition holds:
\begin{enumerate}
	\item[$(A_1)$] The Hill's equation $z''+a_i(t)z=0$ is non-resonant (i.e., its unique periodic solution is the trivial one) and the Green's function $G_i$ has constant sign, that is, either
	\begin{enumerate}
		\item[$(a)$] $G_i(t,s)>0$ for all $(t,s)\in[0,T]\times [0,T]$; or
		\item[$(b)$] $G_i(t,s)<0$ for all $(t,s)\in[0,T]\times [0,T]$. 
	\end{enumerate}
\end{enumerate}

Note that if $a_i$ is not constant there is a well-known $L^p$-criterion which guarantees the constant sign of the Green's function established in \cite{Torres}. Further improvements and complementary results in this line can also be found in \cite{CC_London,CCL}. In the constant case $a_i(t)\equiv k^2$, with $0<k< \pi/T$, the explicit expression of the Green's function is well-known and it is positive, see for instance \cite{Cabada}.

Clearly, system \eqref{eq_periodic_sys} is equivalent to the following system of Hammerstein type integral equations
\begin{align*}
	x(t)=\int_{0}^{T}G_1(t,s)f_1(s,x(s),y(s))\,ds=:\mathcal{T}_1(x,y)(t), \quad t\in[0,T], \\ 
	y(t)=\int_{0}^{T}G_2(t,s)f_2(s,x(s),y(s))\,ds=:\mathcal{T}_2(x,y)(t), \quad t\in[0,T],
\end{align*}
where $G_i$ denote the corresponding Green's function associated to $z''+a_i(t)z=0$, $i=1,2$. Then we look for fixed points of the operator $\mathcal{T}=(\mathcal{T}_1,\mathcal{T}_2):K_1\times K_2\rightarrow X\times X$, where $X$ is the Banach space of continuous functions $\mathcal{C}[0,T]$ with the usual sup-norm $\left\|u\right\|_{\infty}:=\max_{t\in[0,T]}\left|u(t)\right|$ and 
\[K_i:=\left\{u\in X:\min_{t\in[0,T]}u(t)\geq c_i\left\|u\right\|_{\infty} \right\} \qquad (i=1,2), \]
with $c_i:=\min\{m_i/M_i, M_i/m_i \}$, $m_i:=\min_{t,s\in[0,T]}G_i(t,s)$ and $M_i:=\max_{t,s\in[0,T]}G_i(t,s)$. Observe that if $G_i$ is positive, then $0<m_i<M_i$ and so $c_i=m_i/M_i\in (0,1)$, whereas if $G_i$ is negative, then $m_i<M_i<0$ and $c_i=M_i/m_i\in (0,1)$.

Since we are interested in positive solutions for \eqref{eq_periodic_sys}, let us assume that for each $i\in\{1,2\}$ the nonlinearity $f_i$ fulfills the following sign condition:
\begin{enumerate}
	\item[$(A_2)$] The function $f_i$ is non-negative if $(A_1)$-$(a)$ holds and it is non-positive in case that condition $(A_1)$-$(b)$ is satisfied.
\end{enumerate}
Note that, under assumptions $(A_1)$ and $(A_2)$, it is standard to check that the operator $\mathcal{T}$ maps the cone $K:=K_1\times K_2$ into itself, see \cite{Torres} for a similar reasoning. Moreover, by the continuity of the Green's functions and the nonlinearities, it follows that $\mathcal{T}$ is completely continuous (i.e., it is continuous and maps bounded sets into relatively compact ones) as a straightforward consequence of Ascoli-Arzel\'a theorem.

\begin{theorem}\label{th_exist_per}
	Assume that for each $i\in\{1,2\}$ conditions $(A_1)$ and $(A_2)$ hold so that the Green's funtion $G_i$ is positive and, moreover, there exist $\alpha_i,\beta_i>0$ with $\alpha_i\neq \beta_i$ such that the following conditions are fulfilled
	\begin{align}\label{cond_alpha}
		& f_i(t,x_1,x_2)\geq \dfrac{1}{m_i T c_i}x_i \quad \text{for all } x_i\in[c_i \alpha_i,\alpha_i], \ x_j\in[c_j r_j,R_j] \ (j\neq i), \ t\in[0,T]; \text{ and} \\
		& f_i(t,x_1,x_2)\leq \dfrac{1}{M_i T}x_i \quad \text{for all } x_i\in[c_i \beta_i,\beta_i], \ x_j\in[c_j r_j,R_j] \ (j\neq i), \ t\in[0,T], \label{cond_beta}
	\end{align}
	where $r_i:=\min\{\alpha_i,\beta_i \}$ and $R_i:=\max\{\alpha_i,\beta_i \}$, $i=1,2$.
	
	Then problem \eqref{eq_periodic_sys} has at least one positive solution $(x,y)$ with $r_1\leq\left\|x\right\|_{\infty}\leq R_1$ and $r_2\leq\left\|y\right\|_{\infty}\leq R_2$.
\end{theorem}

\noindent
{\bf Proof.} Consider the operator $\mathcal{T}=(\mathcal{T}_1,\mathcal{T}_2):\overline{K}_{r,R}\rightarrow K$ defined as
\[\mathcal{T}_i(x_1,x_2)(t):=\int_{0}^{T}G_i(t,s)f_i(s,x_1(s),x_2(s))\,ds, \quad t\in [0,T] \quad (i=1,2). \]
Let us see that $\mathcal{T}$ is under the assumptions of Theorem \ref{th} and so it has a fixed point $(x,y)\in \overline{K}_{r,R}$, which is a solution to the system \eqref{eq_periodic_sys}.

First, for each $i\in\{1,2\}$, let us prove that $\left\|\mathcal{T}_i(x_1,x_2)\right\|_{\infty}\geq \left\|x_i\right\|_{\infty}$ if $\left\|x_i\right\|_{\infty}=\alpha_i$. Indeed, if $(x_1,x_2)\in\overline{K}_{r,R}$ with $\left\|x_i\right\|_{\infty}=\alpha_i$, then $c_i\alpha_i\leq x_i(t)\leq \alpha_i$ and $c_j r_j\leq x_j(t)\leq R_j$ ($j\neq i$) for all $t\in [0,T]$. Hence, by \eqref{cond_alpha}, we have for each $t\in[0,T]$ that
	\[\mathcal{T}_i(x_1,x_2)(t)=\int_{0}^{T}G_i(t,s)f_i(s,x_1(s),x_2(s))\,ds
	\geq m_i\int_{0}^{T}f_i(s,x_1(s),x_2(s))\,ds \geq \dfrac{1}{T c_i}\int_{0}^{T}x_i(s)\,ds\geq \alpha_i,\]
which clearly implies that $\left\|\mathcal{T}_i(x_1,x_2)\right\|_{\infty}\geq \left\|x_i\right\|_{\infty}$ provided that $\left\|x_i\right\|_{\infty}=\alpha_i$.

Second, for every $i\in\{1,2\}$, let us show that $\left\|\mathcal{T}_i(x_1,x_2)\right\|_{\infty}\leq \left\|x_i\right\|_{\infty}$ if $\left\|x_i\right\|_{\infty}=\beta_i$. Note that, if $(x_1,x_2)\in \overline{K}_{r,R}$ with $\left\|x_i\right\|_{\infty}=\beta_i$, then $c_i\beta_i\leq x_i(t)\leq \beta_i$ and $c_j r_j\leq x_j(t)\leq R_j$ ($j\neq i$) for all $t\in [0,T]$. Thus, we deduce from \eqref{cond_beta} that for every $t\in [0,T]$,
\[\mathcal{T}_i(x_1,x_2)(t)=\int_{0}^{T}G_i(t,s)f_i(s,x_1(s),x_2(s))\,ds\leq M_i\int_{0}^{T}f_i(s,x_1(s),x_2(s))\,ds\leq \dfrac{1}{T}\int_{0}^{T}x_i(s)\,ds\leq \beta_i, \]
and so the inequality $\left\|\mathcal{T}_i(x_1,x_2)\right\|_{\infty}\leq \left\|x_i\right\|_{\infty}$ follows from the definition of the sup-norm.

Finally, the conclusion is a consequence of Theorem \ref{th}, considering that each operator $\mathcal{T}_i$ is compressive if $\alpha_i<\beta_i$ and it is of expansive type when $\alpha_i>\beta_i$.
\qed

Observe that the existence of the positive numbers $\alpha_i$ and $\beta_i$ ($i=1,2$) satisfying conditions \eqref{cond_alpha} and \eqref{cond_beta} can be ensured provided that the nonlinearities exhibit an appropriate asymptotic behavior at zero and at infinity. In order to present such a result, let us introduce the following notations:
\begin{align*}
	(f_i)_{0}:=\lim_{x_i\to 0^+}\dfrac{f_i(t,x_1,x_2)}{x_i} \quad \text{ and } \quad
	(f_i)_{\infty}:=\lim_{x_i\to +\infty}\dfrac{f_i(t,x_1,x_2)}{x_i}.
\end{align*}

\begin{corollary}\label{cor_exist_lim}
Assume that for each $i\in\{1,2\}$ conditions $(A_1)$ and $(A_2)$ hold so that the Green's function $G_i$ is positive and furthermore, one of the following conditions is satisfied:
	\begin{enumerate}
		\item[$(a)$] $(f_i)_{0}=+\infty$ and $(f_i)_{\infty}=0$, uniformly for  $t$, $ x_j$ $ (j\neq i)$; or
		\item[$(b)$] $(f_i)_{0}=0$ and $(f_i)_{\infty}=+\infty$, uniformly for  $t$, $ x_j$ $ (j\neq i)$. 
	\end{enumerate}
	Then problem \eqref{eq_periodic_sys} has at least one positive solution.
\end{corollary}

\noindent
{\bf Proof.} Suppose that condition $(a)$ is satisfied for both components (the proof is similar when condition $(b)$ holds instead of $(a)$). By $(f_i)_{0}=+\infty$ uniformly for $t$, $x_j$ $(j\neq i)$, there exists $\alpha_i>0$ such that
\[\dfrac{f_i(t,x_1,x_2)}{x_i}\geq \dfrac{1}{m_i T c_i} \quad \text{for all } x_i\in (0,\alpha_i], \ x_j\in \mathbb{R}_{+} \ (j\neq i), \ t\in[0,T], \]
which clearly implies \eqref{cond_alpha}. On the other hand, it follows from $(f_i)_{\infty}=0$ uniformly for $t$, $x_j$ $(j\neq i)$ that, there exists $\beta_i>\alpha_i$ large enough such that
\[\dfrac{f_i(t,x_1,x_2)}{x_i}\leq \dfrac{1}{M_i T} \quad \text{for all } x_i\in [c_i \beta_i,+\infty), \ x_j\in \mathbb{R}_{+} \ (j\neq i), \ t\in[0,T], \]
and so condition \eqref{cond_beta} holds. Therefore, the conclusion is derived from Theorem \ref{th_exist_per}.
\qed

Note that, following the usual terminology in the literature, the nonlinearity $f_i$ is said to be \textit{sublinear} if it satisfies condition $(a)$ in the previous result, whereas that it is called \textit{superlinear} when $(b)$ holds. Next example illustrates the case in which $f_1$ is sublinear and $f_2$ is superlinear, which is not covered by classical existence results. 

\begin{example}\label{ex_period_subsup}
	The coupled system of second-order differential equations
		\[\left\{\begin{array}{l} x''+a_1\,x=x^{\lambda}+\sin^2(t+y), \\ y''+a_2\,y=(2+\cos x)\,y^{\mu}, \end{array} \right. \]
	has at least one positive $2\pi$-periodic solution $(x,y)$ provided that $0<a_1,a_2<1/4$, $0\leq\lambda<1$ and $\mu>1$. Indeed, the inequalities $0<a_1,a_2<1/4$ ensure that the corresponding Hill's equations are non-resonant and the Green's functions are positive in $[0,2\pi]\times [0,2\pi]$. Moreover, since $0\leq\lambda<1$ and $\mu>1$, the nonlinearities
	\[f_1(t,x,y)=x^{\lambda}+\sin^2(t+y) \ \text{ and } \ f_2(t,x,y)=(2+\cos x)\,y^{\mu}, \quad (t,x,y)\in [0,2\pi]\times\mathbb{R}_{+}^2, \]
	satisfy that
	\[(f_1)_{0}=+\infty, \ (f_1)_{\infty}=0 \text{ uniformly for } t,\, x_2 \quad \text{and} \quad (f_2)_{0}=0, \ (f_2)_{\infty}=+\infty \text{ uniformly for } t,\, x_1. \]
	Consequently, Corollary \ref{cor_exist_lim} ensures that the system has at least one positive $2\pi$-periodic solution.     
\end{example}

We highlight that the previous existence criteria, established for positive Green's functions, work in a similar way if they are both negative or if one of them is positive and the other one is negative.

\begin{theorem}
	Assume that for each $i\in\{1,2\}$ conditions $(A_1)$ and $(A_2)$ hold so that the Green's funtion $G_i$ is negative and, moreover, there exist $\alpha_i,\beta_i>0$ with $\alpha_i\neq \beta_i$ such that the following conditions are fulfilled
	\begin{align}\label{cond_alpha_neg}
		& f_i(t,x_1,x_2)\leq \dfrac{1}{M_i T c_i}x_i \quad \text{for all } x_i\in[c_i \alpha_i,\alpha_i], \ x_j\in[c_j r_j,R_j] \ (j\neq i), \ t\in[0,T]; \text{ and} \\
		& f_i(t,x_1,x_2)\geq \dfrac{1}{m_i T}x_i \quad \text{for all } x_i\in[c_i \beta_i,\beta_i], \ x_j\in[c_j r_j,R_j] \ (j\neq i), \ t\in[0,T], \label{cond_beta_neg}
	\end{align}
	where $r_i:=\min\{\alpha_i,\beta_i \}$ and $R_i:=\max\{\alpha_i,\beta_i \}$, $i=1,2$.
	
	Then problem \eqref{eq_periodic_sys} has at least one positive solution $(x,y)$ with $r_1\leq\left\|x\right\|_{\infty}\leq R_1$ and $r_2\leq\left\|y\right\|_{\infty}\leq R_2$.
\end{theorem} 

\begin{theorem}\label{th_period_2sign}
	Assume that conditions $(A_1)$ and $(A_2)$ hold so that $G_1$ is positive and $G_2$ is negative. If there exist $\alpha_i,\beta_i>0$ with $\alpha_i\neq \beta_i$ ($i=1,2$) such that \eqref{cond_alpha} and \eqref{cond_beta} hold for $i=1$ and \eqref{cond_alpha_neg} and \eqref{cond_beta_neg} hold for $i=2$, then problem \eqref{eq_periodic_sys} has at least one positive solution $(x,y)$ with $r_1\leq\left\|x\right\|_{\infty}\leq R_1$ and $r_2\leq\left\|y\right\|_{\infty}\leq R_2$.
\end{theorem}

\begin{remark}
	Clearly, a sufficient condition for the existence of $\alpha_i$ and $\beta_i$ satisfying \eqref{cond_alpha_neg} and \eqref{cond_beta_neg}, respectively, is given by the following asymptotic behavior of $f_i$ at zero and infinity:
	\begin{enumerate}
		\item[$(a)$] $(f_i)_{0}=-\infty$ and $(f_i)_{\infty}=0$, uniformly for  $t$, $ x_j$ $ (j\neq i)$; or
		\item[$(b)$] $(f_i)_{0}=0$ and $(f_i)_{\infty}=-\infty$, uniformly for  $t$, $ x_j$ $ (j\neq i)$.
	\end{enumerate}
\end{remark}

Notice that, in the previous results, all the assumptions involving the nonlinearities $f_1$ and $f_2$ can be restricted to a set of the form $[0,T]\times [c_1 r_1, R_1]\times [c_2 r_2, R_2]$. This paves the way for considering systems of \textit{singular} differential equations. Our attention is directed towards nonlinearities with singularities at $x=0$ and $y=0$. 

Let us focus on positive $T$-periodic solutions for a system of the form
\begin{equation}\label{eq_sing}
	\left\{\begin{array}{l} x''+a_1(t)\,x=\dfrac{g_1(t,x,y)}{x^{\lambda}}+h_1(t,x,y), \\[6pt] y''+a_2(t)\,y=-\dfrac{g_2(t,x,y)}{y^{\mu}}-h_2(t,x,y), \end{array} \right.
\end{equation}
where $a_1$, $a_2$ are such that assumption $(A_1)$ holds with $G_1>0$ and $G_2<0$, $\lambda,\mu>0$ and $g_1,g_2,h_1,h_2:[0,T]\times \mathbb{R}_{+}^2\to\mathbb{R}_{+}$ are continuous functions.

\begin{theorem}\label{th_sing}
	Assume that the functions $g_1$ and $g_2$ are bounded below and above by two positive constants $A$ and $B$ and, moreover, 
	\[(h_1)_{\infty}=0, \text{ uniformly for }t,\,y\quad \text{ and } \quad  (h_2)_{\infty}=0, \text{ uniformly for }t,\,x.\]
	Then system \eqref{eq_sing} has at least one positive $T$-periodic solution.	
\end{theorem}

\noindent
{\bf Proof.} Consider the system \eqref{eq_periodic_sys} with the nonlinearities
\[f_1(t,x,y)=\dfrac{g_1(t,x,y)}{x^{\lambda}}+h_1(t,x,y), \quad f_2(t,x,y)=-\dfrac{g_2(t,x,y)}{y^{\mu}}-h_2(t,x,y), \quad (t,x,y)\in [0,T]\times (0,+\infty)^2. \]
Note that assumption $(A_2)$ is satisfied since $f_1\geq 0$ and $f_2\leq 0$. Moreover, for every $t\in [0,T]$ and $x,y>0$,
\[f_1(t,x,y)\geq \dfrac{g_1(t,x,y)}{x^{\lambda}}\geq \dfrac{A}{x^{\lambda}} \ \text{ and } \ f_2(t,x,y)\leq -\dfrac{g_2(t,x,y)}{y^{\mu}}\leq -\dfrac{A}{y^{\mu}}, \]
and thus 
\[\lim_{x\to 0^{+}} f_1(t,x,y)=+\infty \quad \text{uniformly for } t,\ y; \quad \lim_{y\to 0^{+}} f_2(t,x,y)=-\infty \quad \text{uniformly for } t,\ x, \]
which clearly imply the existence of positive numbers $\alpha_1$ and $\alpha_2$ sufficiently close to zero satisfying conditions \eqref{cond_alpha} and \eqref{cond_alpha_neg}, respectively.

On the other hand, it follows from the upper bound of $g_1$ and $g_2$ together with the asymptotic behavior at infinity of $h_1$ and $h_2$ that 
\[(f_1)_{\infty}=0, \text{ uniformly for }t,\,y\quad \text{ and } \quad  (f_2)_{\infty}=0, \text{ uniformly for }t,\,x.\]
Hence, there exist $\beta_1,\beta_2>0$ large enough satisfying conditions \eqref{cond_beta} and \eqref{cond_beta_neg}, respectively. Finally, the conclusion is derived from Theorem~\ref{th_period_2sign}. 
\qed

\begin{example}
Consider the singular system 
	\begin{equation*}
		\left\{\begin{array}{l} x''+a_1\,x=\dfrac{1}{(1+e^{-y}) x^{\mu}}+\dfrac{q_2(t)}{y^2+1}x^{\nu}, \\[10pt] y''+a_2\,y=-\dfrac{2+\cos(x)}{y^{\lambda}}-q_1(t), \end{array} \right.
	\end{equation*}
with $a_1<0$, $0<a_2<1/4$, $\lambda,\mu>0$, $0\leq\nu<1$ and non-negative continuous functions $q_1$ and $q_2$. Then, this system has at least one $2\pi$-periodic positive solution, as a direct consequence of Theorem \ref{th_sing}.	
\end{example}

\begin{remark}\label{rmk_a_r}
	Notice that the singular nonlinearities considered here satisfy that
	\[\lim_{x\to 0^{+}}f_1(t,x,y)=+\infty \quad \text{and} \quad \lim_{y\to 0^{+}}f_2(t,x,y)=-\infty. \]	
	Following the standard terminology in the related literature, it is said that $f_1$ has a \textit{repulsive singularity}, whereas $f_2$ has an \textit{attractive singularity}. We emphasize that this hybrid behavior of the nonlinear parts has not been considered in previous works on systems of singular periodic differential equations such as \cite{CTZ,ft,WCS,W}.
\end{remark}

\section*{Acknowledgements}

J. Rodr\'iguez--L\'opez has been partially supported by Ministerio de Ciencia y Tecnología (Spain), AEI and Feder, grant PID2020-113275GB-I00, and Xunta de Galicia, Spain, Project ED431C 2023/12.

\end{document}